\documentclass[a4paper,12pt]{amsart}%
\usepackage[a4paper,left=1.2cm,right=1.2cm,top=1.45cm,bottom=1.45cm,
headheight=0cm,headsep=0.75cm,includehead,includefoot,footskip=0.75cm,
marginparsep=0cm,marginparwidth=0cm]{geometry}
\usepackage{amsaddr}
\usepackage{amsmath}
\usepackage{amsfonts}
\usepackage{amssymb}
\usepackage{graphicx}
\usepackage{subcaption}
\usepackage{tikz}
\usepackage[colorlinks=true]{hyperref}%
\newlength{\tikzdrawunit}
\newlength{\tikzunit}
\def\pointfillcolor{orange}
\def\smallpointfillcolor{black}
\tikzset{%
  axesstyle/.style={line width=\tikzdrawunit,draw=black},
  projectionstyle/.style={line width=2*\tikzdrawunit,draw=gray,dash pattern=on 0.5\tikzunit off 0.25\tikzunit},
  linestyle/.style={line width=3*\tikzdrawunit,draw=gray},
  contourstyle/.style={line width=4*\tikzdrawunit,draw=black},
  gridstyle/.style={line width=\tikzdrawunit,draw=gray},
  clipstyle/.style={rounded corners=\tikzdrawunit},
  shapestyle/.style={rounded corners=\tikzdrawunit,line width=4*\tikzdrawunit,draw=black},
  inside/.style 2 args={postaction={clip, postaction={draw=#2,line width=#1*2}},line width=0pt},
  outside/.style 2 args={draw=#2,line width=2*#1*\tikzdrawunit},
  pointstyle/.style={fill=\pointfillcolor,draw=white,line width=2*\tikzdrawunit,radius=5*\tikzdrawunit},
  smallpointstyle/.style={fill=\smallpointfillcolor,draw=white,line width=\tikzdrawunit,radius=2*\tikzdrawunit},
  x=\tikzunit,y=\tikzunit,step=1}
\newcommand{\roundmark}[2]{\path[pointstyle] (#1,#2) circle;}
\newcommand{\smallroundmark}[2]{\path[smallpointstyle] (#1,#2) circle;}
\newcommand{\upright}[1]{foreach \x in {1,...,#1}{|- ++(1,1)}}

\newcommand{\downright}[1]{foreach \x in {1,...,#1}{|- ++(1,-1)}}
\newcommand{\rightdown}[1]{foreach \x in {1,...,#1}{-| ++(1,-1)}}
\newcommand{\upleft}[1]{foreach \x in {1,...,#1}{|- ++(-1,1)}}
\newcommand{\leftup}[1]{foreach \x in {1,...,#1}{-| ++(-1,1)}}
\newcommand{\downleft}[1]{foreach \x in {1,...,#1}{|- ++(-1,-1)}}

\newcommand{\staircaseUR}[3]{(#1,#2) -- (\numexpr #1 - #3\relax,#2)\upright{#3} -- cycle}
\newcommand{\staircaseUL}[3]{(#1,#2) -- (\numexpr #1 + #3\relax,#2)\upleft{#3} -- cycle}
\newcommand{\staircaseDR}[3]{(#1,#2) -- (\numexpr #1 - #3\relax,#2)\downright{#3} -- cycle}
\newcommand{\staircaseDL}[3]{(#1,#2) -- (\numexpr #1 + #3\relax,#2)\downleft{#3} -- cycle}
\newcommand{\aztecdiamond}[3]{(\numexpr #1 - #3\relax,#2)\upright{#3}\rightdown{#3}\downleft{#3}\leftup{#3}--cycle}
\newcommand{\aztecdiamondhalf}[3]{(\numexpr #1 - #3\relax,#2)\upright{#3}\rightdown{#3}--cycle}
\newcommand{\biscuit}[3]{(\numexpr #1 - #3 +1\relax,#2)\upright{#3}\downright{\numexpr #3 - 1\relax}\downleft{#3}\upleft{\numexpr #3 - 1\relax}--cycle}
\newcommand{\biscuithalf}[3]{(\numexpr #1 - #3 +1\relax,#2)\upright{#3}\downright{\numexpr #3 - 1\relax}|-cycle}
\newcommand{\StaircaseURGrid}[3]{%
    \begin{scope}
        \clip[clipstyle]\staircaseUR{#1}{#2}{#3};
        \draw[gridstyle] (#1,#2) grid (\numexpr #1 - #3,\numexpr #2 + #3\relax);
    \end{scope}%
}
\newcommand{\StaircaseULGrid}[3]{%
    \begin{scope}
        \clip[clipstyle]\staircaseUL{#1}{#2}{#3};
        \draw[gridstyle] (#1,#2) grid (\numexpr #1 + #3,\numexpr #2 + #3\relax);
    \end{scope}%
}
\newcommand{\StaircaseDRGrid}[3]{%
    \begin{scope}
        \clip[clipstyle]\staircaseDR{#1}{#2}{#3};
        \draw[gridstyle] (#1,#2) grid (\numexpr #1 - #3,\numexpr #2 - #3\relax);
    \end{scope}%
}
\newcommand{\StaircaseDLGrid}[3]{%
    \begin{scope}
        \clip[clipstyle]\staircaseDL{#1}{#2}{#3};
        \draw[gridstyle] (#1,#2) grid (\numexpr #1 + #3,\numexpr #2 - #3\relax);
    \end{scope}%
}
\newcommand{\AztecDiamondGrid}[3]{%
    \begin{scope}
        \clip[clipstyle]\aztecdiamond{#1}{#2}{#3};
        \draw[gridstyle] (\numexpr #1 - #3\relax,\numexpr #2 - #3\relax) grid (\numexpr #1 + #3\relax,\numexpr #2 + #3\relax);
    \end{scope}%
}
\newcommand{\AztecDiamondHalfGrid}[3]{%
    \begin{scope}
        \clip[clipstyle]\aztecdiamondhalf{#1}{#2}{#3};
        \draw[gridstyle] (\numexpr #1 - #3\relax,#2) grid (\numexpr #1 + #3\relax,\numexpr #2 + #3\relax);
    \end{scope}%
}
\newcommand{\BiscuitGrid}[3]{%
    \begin{scope}
        \clip[clipstyle]\biscuit{#1}{#2}{#3};
        \draw[gridstyle] (\numexpr #1 - #3\relax,\numexpr #2 - #3 + 1\relax) grid (\numexpr #1 + #3 - 1\relax,\numexpr #2 + #3\relax);
    \end{scope}%
}
\newcommand{\BiscuitHalfGrid}[3]{%
    \begin{scope}
        \clip[clipstyle]\biscuithalf{#1}{#2}{#3};
        \draw[gridstyle] (\numexpr #1 - #3\relax,0) grid (\numexpr #1 + #3 - 1\relax,\numexpr #2 + #3\relax);
    \end{scope}%
}

\newcommand{\StaircaseURwithGrid}[3]{%
    \StaircaseURGrid{#1}{#2}{#3}%
    \draw[shapestyle]\staircaseUR{#1}{#2}{#3};
}
\newcommand{\StaircaseULwithGrid}[3]{%
    \StaircaseULGrid{#1}{#2}{#3}%
    \draw[shapestyle]\staircaseUL{#1}{#2}{#3};
}
\newcommand{\StaircaseDRwithGrid}[3]{%
       \StaircaseDRGrid{#1}{#2}{#3}%
       \draw[shapestyle]\staircaseDR{#1}{#2}{#3};
}
\newcommand{\StaircaseDLwithGrid}[3]{%
    \StaircaseDLGrid{#1}{#2}{#3}%
    \draw[shapestyle]\staircaseDL{#1}{#2}{#3};
}
\newcommand{\AztecDiamondWithGrid}[3]{%
    \AztecDiamondGrid{#1}{#2}{#3}
    \draw[shapestyle]\aztecdiamond{#1}{#2}{#3};
}
\newcommand{\AztecDiamondHalfWithGrid}[3]{%
    \AztecDiamondHalfGrid{#1}{#2}{#3}
    \draw[shapestyle]\aztecdiamondhalf{#1}{#2}{#3};
}
\newcommand{\BiscuitWithGrid}[3]{%
    \BiscuitGrid{#1}{#2}{#3}
    \draw[shapestyle]\biscuit{#1}{#2}{#3};

}
\newcommand{\BiscuitHalfWithGrid}[3]{%
    \BiscuitHalfGrid{#1}{#2}{#3}
    \draw[shapestyle]\biscuithalf{#1}{#2}{#3};

}

\begin{document}
\title{Counting the lattice rectangles inside Aztec diamonds and square biscuits}
\author{Teofil Bogdan}
\address{``Ioan Bob'' Secondary School, Cluj-Napoca, Romania}
\email{teofilbogdan1@gmail.com}
\author{Mircea Dan Rus}
\address{Technical University of Cluj-Napoca, Romania}
\email{rus.mircea@math.utcluj.ro}
\date{\today}

\begin{abstract}
We are counting the lattice rectangles that can be constructed inside several
planar shapes and identify the corresponding sequences in the OEIS.

\end{abstract}
\maketitle

\section{Introduction}

A point $(x,y)\in\mathbb{Z}\times\mathbb{Z}$ is called a \emph{lattice point}.
A set $[x,x^{\prime}]\times\lbrack y,y^{\prime}]\subseteq\mathbb{R}^{2}$, with
$x<x^{\prime}$, $y<y^{\prime}$ and $x,x^{\prime},y,y^{\prime}\in\mathbb{Z}$ is
called a \emph{lattice rectangle}. In particular, a lattice rectangle having
all sides of length $1$ is called a \emph{unit lattice square}.

Let $n$ be a positive integer. An \emph{Aztec diamond of order} $n$ \cite[p.
277]{stanley1999} is obtained by stacking $2n$ rows of consecutive unit
lattice squares, with the centers of rows vertically aligned and consisting
successively of $2,4,\ldots,2n,2n,\ldots,4,2$ squares. This shape is symmetric
with respect to some lattice point $(p,q)$ (that will be called \emph{the
center} of the Aztec diamond) and can also be described as the union of those
unit lattice squares that lie inside the tilted square $\left\{
(x,y)\in\mathbb{R}^{2}:\left\vert x-p\right\vert +\left\vert y-q\right\vert
\leq n+1\right\}  $.

\begin{figure}[hbt]
    \def\m{4}
    \begin{tikzpicture}
        \foreach \n in {1,...,\m}{%
        \edef\p{\numexpr \n+\n*\n}
        \AztecDiamondWithGrid{\p}{0}{\n}
        \roundmark{\p}{0}%
        }%
    \end{tikzpicture}
    \caption{Aztec diamonds of order up to $\m$ and their corresponding centers marked in \pointfillcolor.}
    \label{intro:first-four-aztecs}
\end{figure}
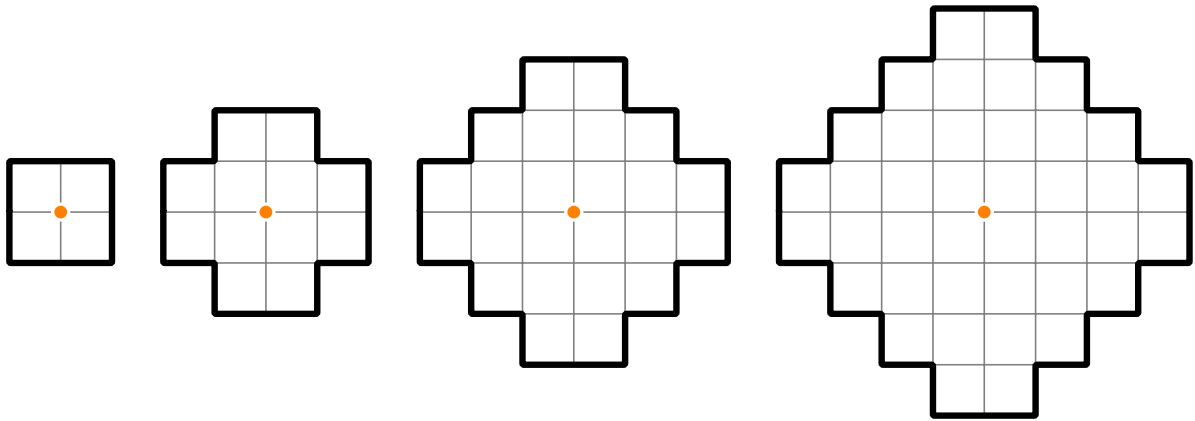 

The \emph{square biscuit of order} $n$ is defined in a similar fashion, by
stacking $2n-1$ rows with their centers vertically aligned which consist
successively of $1,3,\ldots,2n-3,2n-1,2n-3,\ldots,2,1$ consecutive unit
lattice squares. The coordinates of its center are some half integers
$p+\frac{1}{2}$ and $q+\frac{1}{2}$, respectively, so the square biscuit is
the union of the unit lattice squares inside the tilted square $\left\{
(x,y)\in\mathbb{R}^{2}:\left\vert x-p-\frac{1}{2}\right\vert +\left\vert
y-q-\frac{1}{2}\right\vert \leq n\right\}  $. The lattice point $(p,q)$ will
be called \emph{the quasi-center} of the square biscuit.

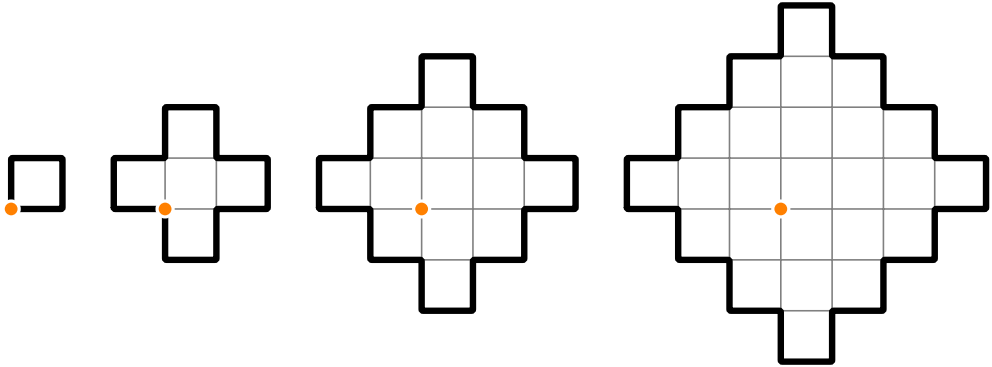
\begin{figure}[h]
    \def\m{4}
    \begin{tikzpicture}
        \draw[shapestyle] (0,0) rectangle (1,1);
        \roundmark{0}{0}
        \foreach \n in {2,...,\m}{%
        \edef\p{\numexpr \n*\n-1}
        \BiscuitWithGrid{\p}{0}{\n}
        \roundmark{\p}{0}%
        }%
    \end{tikzpicture}
    \caption{Biscuits of order up to $\m$ and their corresponding quasi-centers marked in \pointfillcolor.}
    \label{intro:first-four-biscuits}
\end{figure} 

A \emph{staircase of order} $n$ is obtained by stacking $n$ rows of
consecutive unit lattice squares, aligned either to the left or to the right,
which consist of $1,2,3,\ldots,n$ squares and which are stacked either in the
increasing or in the decreasing order of their lengths. Although there are
four\footnote{for $n=1$, there is only one type} types of \emph{staircases},
depending on the alignment of the rows and the ordering in the stack, they are
identical up to a rotation.

\begin{figure}[hb]
    \def\n{4}%
    \hfill
    \begin{subfigure}[b]{0.45\textwidth}
    \centering
    \begin{tikzpicture}
        \StaircaseURwithGrid{0}{0}{\n}
        \StaircaseULwithGrid{0}{0}{\n}
        \StaircaseDLwithGrid{0}{0}{\n}
        \StaircaseDRwithGrid{0}{0}{\n}
        \draw[shapestyle] \aztecdiamond{0}{0}{\n};
    \end{tikzpicture}
    \caption{The four types of staircases of order $\n$. Put together, they form an Aztec diamond of the same order.}
    \label{intro:aztecs2staircases}
    \end{subfigure}
    \hfill%
    \begin{subfigure}[b]{0.45\textwidth}
    \centering
    \begin{tikzpicture}
        \BiscuitWithGrid{0}{0}{\n};
        \draw[contourstyle] (\numexpr 1-\n,0) -- (\n,0);
        \draw[contourstyle] (0,\numexpr 1-\n) -- (0,\n);
        \roundmark{0}{0}%
    \end{tikzpicture}
    \caption{A biscuit of order $\n$ split into four staircases by a vertical line and a horizontal line that intersect at the quasi-center of the biscuit.}
    \label{intro:biscuit2staircases}
    \end{subfigure}
    \caption{Aztec diamonds, biscuits and staircases.}
    \label{intro:staircases}
    \hfill
\end{figure}
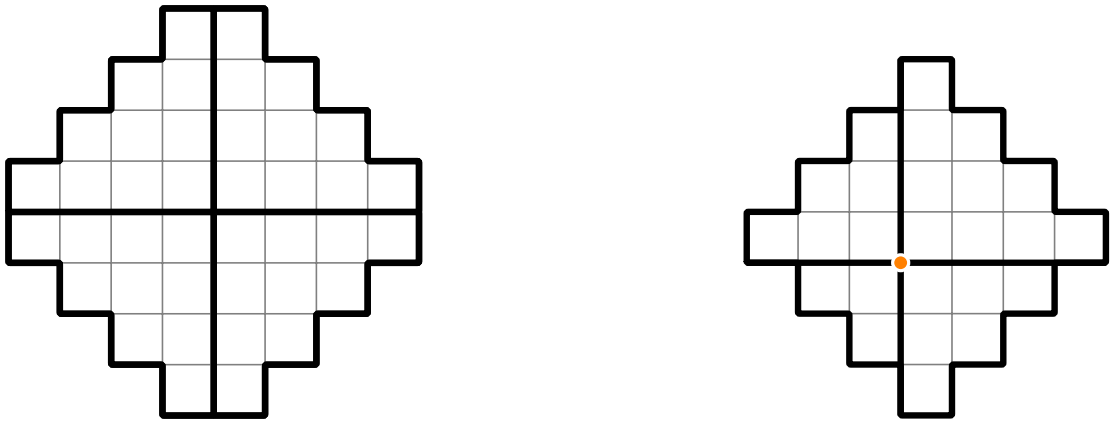

Splitting an Aztec diamond, either vertically or horizontally, through the
center gives two halves that are identical up to a rotation. In the case of a
square biscuit\footnote{of order at least 2}, a vertical or a horizontal
splitting through the quasi-center produces two different halves, one larger
by one row than the other. Going further, both an Aztec diamond and a square
biscuit\footnote{not too small} of order $n$ can be split into four
staircases, one of each type. The staircases in the former case have all the
same order $n$, while for the latter there is one staircase of order $n$, two
of order $n-1$ and one of order $n-2$ (see Figure \ref{intro:staircases}).

\newpage

\section{Problems}

\noindent\textbf{Main problems:} Find the number of lattice
rectangles\footnote{squares included} included in:

\begin{enumerate}
\item[(1)] an Aztec diamond of order $n$;

\item[(2)] a square biscuit of order $n$.
\end{enumerate}

\medskip

\noindent\textbf{Intermediate problems:} Find the number of lattice
rectangles\footnote{again, squares included} included in:

\begin{enumerate}
\item[(3)] a staircase of order $n$;

\item[(4)] a half of an Aztec diamond of order $n$;

\item[(5)] the larger half of a square biscuit of order $n$.
\end{enumerate}

\medskip

We will denote the answers to these five problems by $a(n)$, $b(n)$,
$s(n)$, $a_{\frac{1}{2}}(n)$ and $b_{\frac{1}{2}}(n)$, respectively.

Problems 2, 3 and 5 were studied by the authors in \cite{bogdan2020}. Here, we
present alternative solutions to Problems 2 and 5 that make the connection
with Problems 1 and 4.

\section{Solutions to the intermediate problems}

\begin{proof}
[Solution to Problem 3]Consider a staircase of order $n$, positioned in the
plane as shown in Figure \ref{solutions:staircase}. A lattice rectangle
$[a,b]\times\lbrack c,d]$ is included in the staircase if and only if $0\leq
a<b\leq n$, $1<c<d\leq n+2$ and $b+1\leq c$, where the last condition means
that the bottom right corner of the rectangle does not lie below the line
$y=x+1$. Concluding, any lattice rectangle that lies inside the staircase is
uniquely determined by a quadruple $(a,b,c,d)$ of integers that satisfy $0\leq
a<b<c<d\leq n+2$. There are $\dbinom{n+3}{4}$ such quadruples, hence
$s(n)=\dbinom{n+3}{4}$.

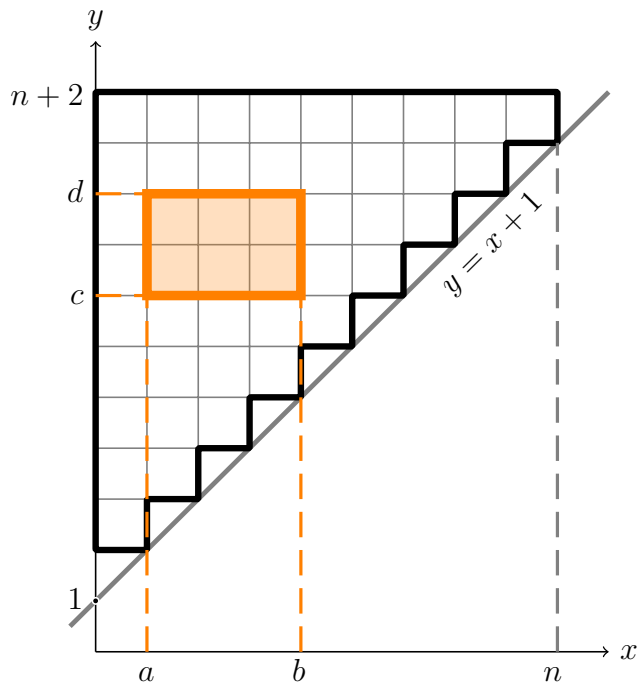
\begin{figure}[ht]
    \def\n{9}\def\a{1}\def\b{4}\def\c{7}\def\d{9}
    \begin{tikzpicture}

        \begin{scope}[axesstyle]
            \draw[->] (0,0) -- (\numexpr 1+\n,0) node[right] {$x$};
            \draw[->] (0,0) -- (0,\numexpr 3+\n) node[above] {$y$};
        \end{scope}
        \draw[linestyle] (-0.5,0.5) -- node[sloped,below,near end] {$y=x+1$} (\numexpr 1+\n,\numexpr 2+\n);
        \StaircaseDLwithGrid{0}{\numexpr 2+ \n}{\n}
        \begin{scope}[projectionstyle,draw=orange]
            \draw (\a,0) node[anchor=north] {\parbox[t][1.25ex][b]{1.25ex}{$a$}} -- (\a,\c);
            \draw (\b,0) node[anchor=north] {\parbox[t][1.25ex][b]{1.25ex}{$b$}} -- (\b,\c);
            \draw (0,\c) node[anchor=east]  {\parbox[t][1.25ex][b]{2ex}{\hfill$c$}} -- (\a,\c);
            \draw (0,\d) node[anchor=east]  {\parbox[t][1.25ex][b]{2ex}{\hfill$d$}} -- (\a,\d);
            \draw[gray] (\n,0) node[anchor=north]  {\parbox[t][1.25ex][b]{2ex}{\color{black} $n$}} -- (\n,\numexpr 1+\n);
            \node[anchor=east] at (0,\numexpr 2+ \n) {\parbox[t][2ex][b]{7ex}{\hfill$n+2$}};
        \end{scope}
        \draw[rounded corners=0\tikzdrawunit,line width=6*\tikzdrawunit,draw=orange,fill=orange,fill opacity=0.25] (\a,\c) rectangle (\b,\d);
        \node[anchor=east] at (0,1) {\parbox[t][1.25ex][b]{1.25ex}{$1$}};
        \smallroundmark{0}{1}
    \end{tikzpicture}
    \caption{Finding the relations between the coordinates of a lattice rectangle inside a staircase of order $n$.}
    \label{solutions:staircase}
\end{figure}

\end{proof}

\begin{proof}
[Solution to Problem 4]Fix $n\geq2$. Denote by $A_{1/2}$ the upper half of an
Aztec diamond of order $n$ and let $\Delta$ be its vertical axis of symmetry
that splits $A_{1/2}$ into two staircases of order $n$.

There are $2s(n)$ lattice rectangles inside the Aztec diamond that lie
entirely either in the left staircase or in the right staircase.

For the remaining lattice rectangles (see Figure
\ref{solutions:half-aztec-diamond-a}), $\Delta$ splits each rectangle into two
smaller ones (that will be individually referred to as \emph{the left part}
and \emph{the right part} of the rectangle). We say that a rectangle is of
type $L$, $R$ or $C$\footnote{these notations are abbreviations of the position
of the rectangle with respect to $\Delta$: left, right or centered}  if its
left part is larger, smaller or equal in size to
its right part, respectively.

\begin{figure}[h]
    \def\n{7}\def\a{-3}\def\b{2}\def\c{1}\def\d{3}
    \begin{tikzpicture}
        \AztecDiamondHalfWithGrid{0}{0}{\n}
        \draw[draw=black, line width = 2*\tikzdrawunit,dotted] (0,-1) node[anchor=south west] {$\Delta$} -- (0,\numexpr 1 +\n);
        \draw[inside={6*\tikzdrawunit}{orange},outside={\tikzdrawunit}{orange},fill=orange,fill opacity=0.25] (\a,\c) rectangle (\b,\d);
    \end{tikzpicture}
    \caption{Half of an Aztec diamond and an $L$--type rectangle}
    \label{solutions:half-aztec-diamond-a}
\end{figure}
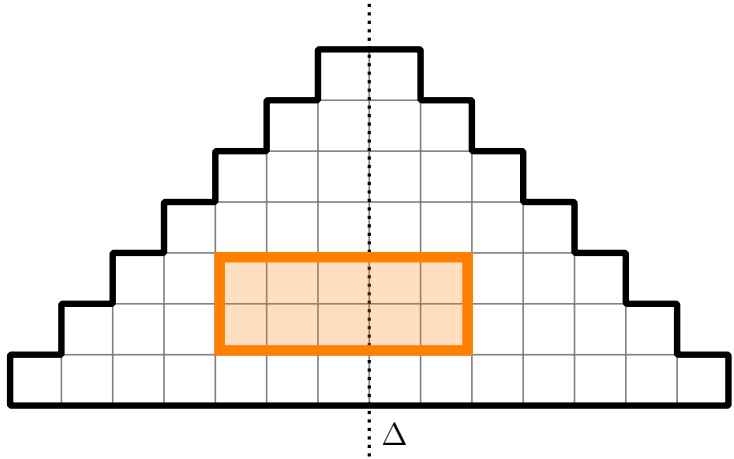 

We count the type--$L$ rectangles by the following bijective argument. Each
such rectangle can be uniquely identified with the difference between the
reflection of its left part with respect to $\Delta$ and its right part (in
Figure \ref{solutions:half-aztec-diamond-b}, the orange type--$L$ rectangle is
transformed into the blue-filled rectangle). The result of this transformation
is always a rectangle that is included in the staircase of order $n-1$,
obtained by the vertical split of $A_{1/2}$ one unit to the right of $\Delta$
and taking the right part. It is easy to check that the transformation is
uniquely reversible, hence bijective, which gives the number of type--$L$
rectangles to be $s(n-1)$.

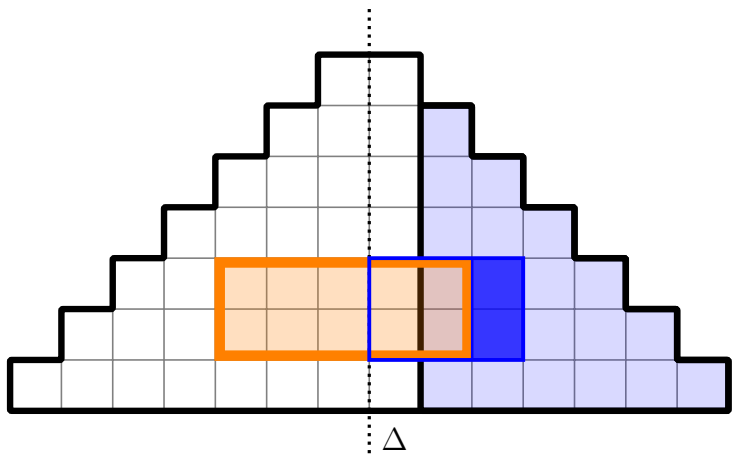
\begin{figure}[h]
    \def\n{7}\def\a{-3}\def\b{2}\def\c{1}\def\d{3}
    \begin{tikzpicture}
        \AztecDiamondHalfWithGrid{0}{0}{\n}
        \draw[draw=black, line width = 2*\tikzdrawunit,dotted] (0,-1) node[anchor=south west] {$\Delta$} -- (0,\numexpr 1 +\n);
        \draw[shapestyle,fill=blue,fill opacity=0.15] \staircaseUL{1}{0}{\numexpr \n -1\relax};
        \draw[inside={6*\tikzdrawunit}{orange},outside={\tikzdrawunit}{orange},fill=orange,fill opacity=0.25] (\a,\c) rectangle (\b,\d);
        \draw[inside={0\tikzdrawunit}{blue},outside={2\tikzdrawunit}{blue}] (0,\c) rectangle (\numexpr -\a,\d);
        \draw[inside={\tikzdrawunit}{blue},outside={\tikzdrawunit}{blue},fill=blue,fill opacity=0.75] (\b,\c) rectangle (\numexpr -\a,\d);
    \end{tikzpicture}
    \caption{A bijective correspondence between the $L$--type rectangles in half of an Aztec diamond of order $n$ and the rectangles in a staircase of order $n-1$. }
    \label{solutions:half-aztec-diamond-b}
\end{figure} 

By symmetry, the number of type--$R$ rectangles is equal to the number of
type--$L$ ones.

Also, every type--$C$ rectangle can be uniquely described by its right part,
which is a rectangle included in the staircase of order $n$ to the right of
$\Delta$, having the left side on $\Delta$. The number of such rectangles is
$s(n)-s(n-1)$, since there are $s(n-1)$ rectangles inside the right staircase
of order $n$ that do not have the left side on $\Delta$.

Concluding, there are $2s(n-1)+s(n)-s(n-1)=s(n)+s(n-1)$ lattice rectangles
included in $A_{1/2}$ whose interior is intersected by $\Delta$, so
\[
a_{\frac{1}{2}}(n)=2s(n)+s(n)+s(n-1)=3s(n)+s(n-1)=3\dbinom{n+3}{4}%
+\dbinom{n+2}{4}=\dfrac{n\left(  n+1\right)  \left(  n+2\right)  ^{2}}{6}%
\]
for all $n\geq2$ and with $a_{\frac{1}{2}}(1)=3$ also satisfying the general
formula (consider that $s(0):=0$).
\end{proof}

\begin{proof}
[Solution to Problem 5]Fix $n\geq2$ . Denote by $B_{1/2}$ the larger upper
half of a square biscuit of order $n$ and by $\delta$ its vertical
axis\footnote{this is not a lattice line} of symmetry. There are $2s(n-1)$
lattice rectangles included in $B_{1/2}$ that lie entirely either to the left
or to the right of $\delta$.

When counting the remaining lattice rectangles whose interior is intersected
by $\delta$ (see Figure \ref{solutions:half-biscuit-a}), it is enough to
expand $B_{1/2}$ to half of an Aztec diamond (denote it by $A_{1/2}$), by
inserting a middle column of $n$ unit lattice squares to the left\footnote{or
to the right} of $\delta$ (see Figure \ref{solutions:half-biscuit-b}).
Naturally, the rectangles that intersect $\delta$ will also expand one column
to the left of $\delta$.

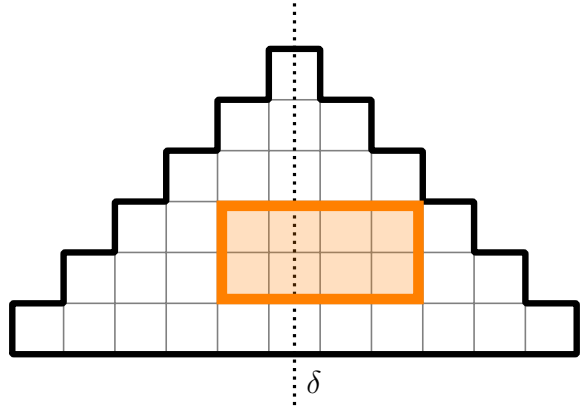
\begin{figure}[ht]
    \def\n{6}\def\a{-1}\def\b{3}\def\c{1}\def\d{3}
    \begin{tikzpicture}
        \BiscuitHalfWithGrid{0}{0}{\n}
        \draw[draw=black, line width = 2*\tikzdrawunit,dotted] (0.5,-1) node[anchor=south west] {$\delta$} -- (0.5,\numexpr 1 +\n);
        \draw[inside={6*\tikzdrawunit}{orange},outside={\tikzdrawunit}{orange},fill=orange,fill opacity=0.25] (\a,\c) rectangle (\b,\d);
    \end{tikzpicture}
    \caption{Half of a square biscuit and a rectangle intersected by $\delta$.}
    \label{solutions:half-biscuit-a}
\end{figure} 

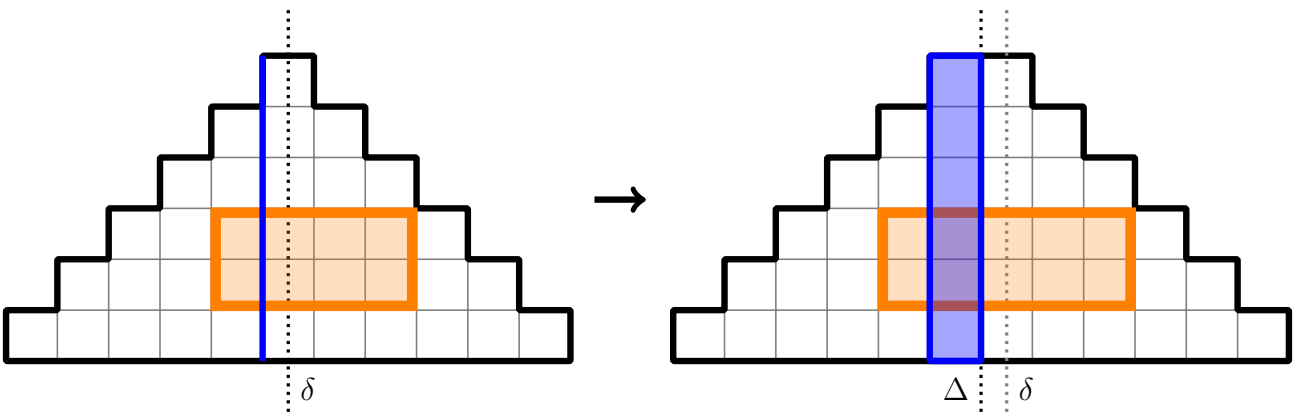
\begin{figure}[ht]
    \def\n{6}\def\a{-1}\def\b{3}\def\c{1}\def\d{3}
    \begin{tikzpicture}
        \BiscuitHalfWithGrid{0}{0}{\n}
        \draw[draw=black, line width = 2*\tikzdrawunit,dotted] (0.5,-1) node[anchor=south west] {$\delta$} -- (0.5,\numexpr 1 +\n);
        \draw[inside={6*\tikzdrawunit}{orange},outside={\tikzdrawunit}{orange},fill=orange,fill opacity=0.25] (\a,\c) rectangle (\b,\d);
        \draw[shapestyle,draw=blue] (0,0) -- (0,\n);
    \end{tikzpicture}
    \begin{tikzpicture}
    \def\p{\numexpr\n/2+1\relax}
    \draw[->,shapestyle] (0,\p) -- (1,\p);
    \node at (0,0) {};
    \end{tikzpicture}
    \def\n{6}\def\a{-2}\def\b{3}\def\c{1}\def\d{3}
    \begin{tikzpicture}
        \AztecDiamondHalfWithGrid{0}{0}{\n}
        \draw[draw=gray, line width = 2*\tikzdrawunit,dotted] (0.5,-1) node[anchor=south west] {$\delta$} -- (0.5,\numexpr 1 +\n);
        \draw[draw=black, line width = 2*\tikzdrawunit,dotted] (0,-1) node[anchor=south east] {$\Delta$} -- (0,\numexpr 1 +\n);
        \draw[inside={6*\tikzdrawunit}{orange},outside={\tikzdrawunit}{orange},fill=orange,fill opacity=0.25] (\a,\c) rectangle (\b,\d);
        \draw[shapestyle,draw=blue,fill=blue,fill opacity=0.35] (-1,0) rectangle (0,\n);
    \end{tikzpicture}
    \caption{Expanding half of a square biscuit to half of an Aztec diamond of the same order, by inserting a middle column of height $n$.}
    \label{solutions:half-biscuit-b}
\end{figure} 

If $\Delta$ denotes the vertical axis of symmetry of
$A_{1/2}$, it is straightforward to check that there exists a bijective
correspondence between the lattice rectangles included in $B_{1/2}$ that
intersect $\delta$ and the lattice rectangles included in $A_{1/2}$ whose
interior is intersected by $\Delta$; their number $s(n)+s(n-1)$ was found in
the solution of Problem 4.

Concluding,%
\[
b_{\frac{1}{2}}(n)=2s(n-1)+s(n)+s(n-1)=s(n)+3s(n-1)=\dbinom{n+3}{4}%
+3\dbinom{n+2}{4}=\frac{n^{2}\left(  n+1\right)  \left(  n+2\right)  }{6}%
\]
for all $n\geq2$ and with $b_{\frac{1}{2}}(1)=1$ also satisfying the general formula.
\end{proof}

\section{Solutions to the main problems}

\begin{proof}
[Solution to Problem 1]Fix $n\geq2$. Let $A$ be an Aztec diamond of order $n$
and $\Delta$ be its vertical axis of symmetry which splits $A$ into two equal
halves. There are $2a_{\frac{1}{2}}(n)$ lattice rectangles inside the Aztec
diamond that lie entirely either in the left half or in the right half.

The remaining lattice rectangles, split in two by $\Delta$, are counted using
the same approach as in the solution of Problem 4, with the staircases
replaced by the halves of $A$, hence $s(n)$ replaced by $a_{\frac{1}{2}}(n)$.
We leave to the interested reader to verify this claim, hence to obtain the
formula $a_{\frac{1}{2}}(n)+a_{\frac{1}{2}}(n-1)$ for the number of the
lattice rectangles included in $A$ and whose interior is intersected by
$\Delta$.

Concluding,
\begin{align*}
a(n)  &  =2a_{\frac{1}{2}}(n)+a_{\frac{1}{2}}(n)+a_{\frac{1}{2}}%
(n-1)=3a_{\frac{1}{2}}(n)+a_{\frac{1}{2}}(n-1)\\
&  =9\dbinom{n+3}{4}+6\dbinom{n+2}{4}+\dbinom{n+1}{4}=\dfrac{n\left(
n+1\right)  \left(  4n^{2}+12n+11\right)  }{6}%
\end{align*}
for all $n\geq2$. A direct count gives $a(1)=9$, which also agrees with the
general formula.
\end{proof}

\begin{proof}
[Solution to Problem 2]Fix $n\geq2$ . Denote by $B$ a square biscuit of order
$n$ and by $\delta$ its vertical axis\footnote{the same line as in the
solution of Problem 5; it is not a lattice line} of symmetry. There are
$2b_{\frac{1}{2}}(n-1)$ lattice rectangles included in $B$ that lie entirely
either to the left or to the right of $\delta$. We leave to the interested
reader to check that the arguments presented in the previous proofs can be
easily adapted in counting the lattice rectangles included in $B$ and
intersected by $\delta$, leading to $b_{\frac{1}{2}}(n)+b_{\frac{1}{2}}(n-1)$ rectangles.

Concluding,%
\begin{align*}
b(n) &  =2b_{\frac{1}{2}}(n-1)+b_{\frac{1}{2}}(n)+b_{\frac{1}{2}%
}(n-1)=b_{\frac{1}{2}}(n)+3b_{\frac{1}{2}}(n-1)\\
&  =\dbinom{n+3}{4}+6\dbinom{n+2}{4}+9\dbinom{n+1}{4}=\dfrac{n\left(
n+1\right)  \left(  4n^{2}-4n+3\right)  }{6}%
\end{align*}
for all $n\geq2$. Also, $b(1)=1$ which agrees with the general formula.
\end{proof}

\section{Identifying the results in the On-Line Encyclopedia of Integer
Sequences}

The sequences $a_{\frac{1}{2}}$, $b_{\frac{1}{2}}$, $a$ and $b$ appear in the
On-Line Encyclopedia of Integer Sequences (OEIS) \cite{oeis} as A004320,
A002417, A330805 and A213840, respectively. At the time this manuscript was
typeset, there was no mention in the OEIS about the combinatorial problems
studied in this paper in connection to the sequences A004320 and A213840. The
OEIS connects the sequences A002417 and A330805 to the corresponding problems
studied in this paper (though using different terminology), but provides no
reference to a proof.

\end{document}